\overfullrule=0pt
\centerline {\bf A minimax theorem in infinite-dimensional topological vector spaces}\par
\bigskip
\bigskip
\centerline {BIAGIO RICCERI}\par
\bigskip
\bigskip
{\bf Abstract:} In this paper, we obtain a minimax theorem by means of which, in turn, we prove the following result:\par
\smallskip
Let $E$ be an infinite-dimensional reflexive real Banach space, $T:E\to E$ a non-zero compact linear operator, 
$\varphi:E\to {\bf R}$ a lower semicontinuous, convex and 
coercive functional, $I\subset {\bf R}$ a compact interval, with $0\in I$, $\psi:I\to {\bf R}$ a lower semicontinuous convex function.\par
Then, for each $r>\varphi(0)$, one has
$$\sup_{x\in X}\inf_{\lambda\in I}(\varphi(T(x)-\lambda x)+\psi(\lambda))=r+\psi(0)\ ,$$
where
$$X=\{x\in E : \varphi(T(x))\leq r\}\ .$$
\bigskip
\bigskip
{\bf Key words:} Minimax theorem; inf-compactness; lower semicontinuity; connectedness; sequential weak topology; compact linear operator.\par
\bigskip
\bigskip
{\bf 2010 Mathematics Subject Classifications:} 49J35; 49K35; 46A55; 46B10; 47A10; 52A41.
\bigskip
\bigskip
\bigskip
\bigskip
Let $E$ be a topological space and $X$ a non-empty subset of $E$. A function $f:X\to {\bf R}$ is said to be relatively inf-compact (resp.
relatively sequentially inf-compact) in $E$,  
provided that, for each
$r\in {\bf R}$, the sub-level set $f^{-1}(]-\infty,r])$ is relatively compact (resp. sequentially relatively compact) in $E$, that is its closure in $E$ 
is compact (resp. sequentially compact).
A real-valued function $f$ on a convex subset of a vector space is said to be quasi-convex if, for each $r\in {\bf R}$, the set $f^{-1}(]-\infty,r])$ is
convex.\par
\smallskip
The aim of this very short note is to highlight the following minimax result:\par
\medskip
THEOREM 1. - {\it Let $E$ be a real Hausdorff topological vector space and let $X\subseteq E$ be an infinite-dimensional convex set whose interior in
its closed affine hull is non-empty. Moreover, let $I\subset {\bf R}$ be a compact interval and $f:X\times I\to {\bf R}$ a function which is
lower semicontinuous in $X\times I$ and quasi-convex in $I$. Finally, assume that there is a set $D\subset I$, dense in $I$, such that, for each
$\lambda\in D$, the function $f(\cdot,\lambda)$ is 
relatively inf-compact (resp. relatively sequentially inf-compact) in $E$\ .\par
Then, one has
$$\sup_{x\in X}\inf_{\lambda\in I}f(x,\lambda)=\inf_{\lambda\in I}\sup_{x\in X}f(x,\lambda)\ .$$}\par
\medskip
Theorem 1 can be qualified as unconventional in the sense that, in most of the known minimax theorems, lower semicontinuity and inf-compactness are related 
to the variable with respect to which one takes the {\it inf}, while it is quasi-concavity that one generally assumes with respect to the other variable (see,
for instance, [4]).
\smallskip
It is natural to ask whether the two assumptions made on the convex set $X$ are necessary. We start just presenting two examples related to such a question.\par
\smallskip
The first example concerns the infinite-dimensionality of $X$.\par
\medskip
EXAMPLE 1. - Let $E$ be a finite-dimensional normed space and let $f:E\times [0,1]\to {\bf R}$ be the function defined by
$$f(x,\lambda)=|\|x\|-\lambda(\|x\|^2+1)|$$
for all $(x,\lambda)\in E\times [0,1]$.\par
Of course, $f$ is convex in $[0,1]$, inf-compact in $E$ and, just because dim$(E)<\infty$, continuous in $E\times [0,1]$. Further, notice that, for each
$x\in E$, taking $\lambda={{\|x\|}\over {\|x\|^2+1}}$, we have $\lambda\in [0,1]$ and $f(x,\lambda)=0$. This implies that
$$\sup_{x\in E}\inf_{\lambda\in [0,1]}f(x,\lambda)=0\ .$$
On the other hand, we clearly have
$$\inf_{\lambda\in [0,1]}\sup_{x\in E}f(x,\lambda)=+\infty\ .$$
So, the conclusion of Theorem 1 can fail if $X$ if finite-dimensional.\par
\medskip
The second example deals with the non-emptyness of the interior of $X$ in its closed affine hull.\par
\medskip
EXAMPLE 2. -  Let $E$ be an infinite-dimensional reflexive real Banach space, let $X$ be the open unit ball in $E$ and let $\varphi\in E^*$, with $\|\varphi\|_{E^*}=1$.
Consider the function $f:X\times [0,1]\to {\bf R}$ defined by
$$f(x,\lambda)=\left | {{1}\over {1-\varphi(x)}}-\lambda\left ( \left ( {{1}\over {1-\varphi(x)}}\right )^2+1\right )\right |$$
for all $(x,\lambda)\in X\times [0,1]$.\par
Consider $E$ equipped with the weak topology. Clearly, the affine hull of $X$ is the whole $E$ and, since dim$(E)=\infty$, the interior of $X$ in the weak topology is empty.
Since, by reflexivity, $X$ is relatively weakly compact, the function $f$ is relatively weakly inf-compact in $E$. Since $\varphi\in E^*$, $f$ is weakly continuous in $X\times
[0,1]$, besides being convex in $[0,1]$. As in Example 1, it is seen that
$$\sup_{x\in X}\inf_{\lambda\in [0,1]}f(x,\lambda)=0$$
and
$$\inf_{\lambda\in [0,1]}\sup_{x\in X}f(x,\lambda)=+\infty\ .$$
So, the conclusion of Theorem 1 can fail if the interior of $X$ in its closed affine hull is empty.\par
\medskip
Our proof of Theorem 1 is fully based on the joint use of three previous results of ours. We now recall them.\par
\medskip
THEOREM A ([2], Proposition 3). - {\it Let $E$ be a real Hausdorff topological vector space, let $X\subseteq E$ be an infinite-dimensional convex set whose interior in
its closed affine hull is non-empty and let  $K\subseteq E$ be a relatively compact  (resp. relatively sequentially compact) set.\par
 Then, the set $X\setminus K$  is connected.}\par
\medskip
REMARK 1. - Notice that, in [2], such a result was proved for the relatively compact case only. The same proof shows the validity of the result also in the
relatively sequentially compact case, in view of the fact that any Hausdorff topological vector space possessing a sequentially compact neighbourhood of
$0$ is finite-dimensional.\par
\medskip
THEOREM B ([3], Proposition 5.3) - {\it Let $X, Y$ be two topological spaces, with $X$
connected,
and let $F:X\to 2^Y$ be a lower semicontinuous multifunction with
non-empty values.
Assume the set
$$\{x\in X : F(x)\hskip 3pt is\hskip 3pt connected\}$$
is dense in $X$.\par
Then, the set
$$\{(x,y)\in X\times Y : y\in F(x)\}$$
is connected.}\par
\medskip
For a generic set $S\subseteq X\times I$, for each $(x,\lambda)\in X\times I$, we set
$$S_x=\{\mu\in I: (x,\mu)\in S\}$$
and
$$S^{\lambda}=\{u\in X : (u,\lambda)\in S\}\ .$$
\medskip
 THEOREM C ([1], Theorem 2.3). - {\it Let $X$ be
 a topological space, $I\subseteq {\bf R}$ a compact interval and $S,T\subseteq
X\times I$. Assume that $S$ is connected and $S^{\lambda}\neq \emptyset$
for all $\lambda\in I$, while $T$ is closed and  $T_x$ is non-empty and connected for
all $x\in X$.\par
Then, one has $S\cap T\neq \emptyset$.}\par
\medskip
{\it Proof of Theorem 1}.
Arguing by contradiction, assume that 
$$\sup_{x\in X}\inf_{\lambda\in I}f(x,\lambda)<\inf_{\lambda\in I}\sup_{x\in X}f(x,\lambda)\ .$$
Fix $\rho$ satisfying
$$\sup_{x\in X}\inf_{\lambda\in I}f(x,\lambda)<\rho<\inf_{\lambda\in I}\sup_{x\in X}f(x,\lambda) \eqno{(1)}$$
and put
$$S=\{(x,\lambda)\in X\times I : f(x,\lambda)>\rho\}$$
and
$$T=\{(x,\lambda)\in X\times I : f(x,\lambda)\leq \rho\}\ .$$
Since $f$ is lower semicontinuous, the set $T$ is closed. Moreover, for each $x\in X$, the set $T_x$ is non-empty by $(1)$ and connected by
the quasi-convexity of $f(x,\cdot)$.  By $(1)$ again, $S^{\lambda}\neq\emptyset$ for all $\lambda\in I$. Fix $\lambda\in D$.  Since
$$S^{\lambda}=X\setminus \{x\in X : f(x,\lambda)\leq\rho\}$$
and $\{x\in X : f(x,\lambda)\leq\rho\}$ is relatively compact (resp. relatively sequentially compact) in $E$, in view of Theorem A, the set 
$S^{\lambda}$ turns out to be connected.
On the other hand, since
$S_x$ is open for all $x\in X$, the multifunction $\lambda\to S^{\lambda}$ is lower semicontinuous in $I$. At this point, we can apply Theorem B
to realize that the set
$$\{(\lambda,x)\in I\times X : (x,\lambda)\in S\}$$
is connected.  But such a set is clearly homeomorphic to $S$, and so $S$ is connected. As a consequence, each assumption of Theorem $C$ is satisfied, and
hence we would have $S\cap T\neq\emptyset$ which is clearly false. Such a contradiction completes the proof.\hfill  $\bigtriangleup$\par
\medskip
We conclude with the following application of Theorem 1. We first introduce a notation. Namely, if $Y$ is a topological space and $\tau$ is the
topology of $Y$, we denote by $\tau_s$ the topology on $Y$ whose members are the sequentially open subsets of $Y$. Let us recall that a set
$A\subseteq Y$ is said to be sequentially open if, for every sequence $\{y_n\}$ in $Y$ converging to a point of $A$, there is $\nu\in {\bf N}$ such that
$y_n\in A$ for all $n\geq\nu$. A functional $\varphi$ on a real normed space is said to be coercive if $\lim_{\|x\|\to +\infty}\varphi(x)=+\infty$.
\medskip
THEOREM 2. - {\it Let $E$ be an infinite-dimensional reflexive real Banach space, $T:E\to E$ a non-zero compact linear operator, 
$\varphi:E\to {\bf R}$ a lower semicontinuous, convex and 
coercive functional, $I\subset {\bf R}$ a compact interval, with $0\in I$, $\psi:I\to {\bf R}$ a lower semicontinuous convex function.\par
Then, for each $r>\varphi(0)$, one has
$$\sup_{x\in X}\inf_{\lambda\in I}(\varphi(T(x)-\lambda x)+\psi(\lambda))=r+\psi(0)\ ,$$
where
$$X=\{x\in E : \varphi(T(x))\leq r\}\ .$$}
\smallskip
PROOF. By a classical result, $\varphi$ turns out to be continuous with respect to the strong topology of $E$. Since $T(E)$ is
a non-zero linear subspace, the set $\varphi(T(E))$ is unbounded above. Indeed, if not, $\varphi$ would be constant on $T(E)$, contrary to the
coercivity of $\varphi$. As a consequence, since $T(E)$ is connected, we have
$$\varphi(T(E))=\left ( \inf_{T(E)}\varphi,+\infty\right [\ .$$
From this, we clearly infer that
$$\sup_{x\in X}\varphi(T(x))=r \ .\eqno{(2)}$$
Next,  consider the function $f:X\times I\to {\bf R}$ defined by
$$f(x,\lambda)=\varphi(T(x)-\lambda x)+\psi(\lambda)$$
for all $(x,\lambda)\in X\times I$. Now, denote by $\tau$ the weak topology of $E$. Notice that $T$, being linear and compact, turns out to be
sequentially continuous from $E$ with the topology $\tau$ to $E$ with the strong topology. It is easy to check that this is equivalent to 
the continuity of $T$ from $E$ with the topology $\tau_s$ to $E$ with the strong topology. Of course, $(E,\tau_s)$ is a Hausdorff topological
vector space. Now, we are going to apply Theorem 1 to the function $f$ considering $E$ with the topology $\tau_s$. First, notice that the set
$X$ is convex and its interior in $\tau_s$ is non-empty. Actually, $X$ contains the non-empty set $T^{-1}(\varphi^{-1}(]-\infty,r[))$ which is
open in $\tau_s$, by the remarks above. Next, observe that, for each $\lambda\in {\bf R}$,
the function $x\to T(x)-\lambda x$, being continuous and linear, is continuous from the weak to the weak topology, and so, {\it a fortiori},
from the $\tau_s$ to the weak topology. Of course, this implies that the function $(x,\lambda)\to T(x)-\lambda x$ is continuous from the product
of $\tau_s$ and the topology of ${\bf R}$ to the weak topology. But then, since $\varphi$ is weakly lower semicontinuous, the function $f$
is lower semicontinuous in $X\times I$ with respect to the considered topology. Of course, $f$ is convex in $I$. Finally, by a classical result,
the spectrum of $T$ is countable, and so the set, say $D$, of all $\lambda\in I$ such that $x\to T(x)-\lambda x$ is a homeomorphism between $E$ 
(with the strong topology) and itself is dense in $I$. Fix $\lambda\in D$. Of course, since $\varphi$ is coercive, for each $\rho\in {\bf R}$, the set
$$\{x\in E : \varphi(T(x)-\lambda x)\leq \rho\}$$
is bounded. Hence, due to the reflexivity of $E$, the sub-level sets of $f(\cdot,\lambda)$ are 
weakly compact and so, by the Eberlein-Smulyan theorem, 
sequentially weakly compact which is equivalent to sequentially $\tau_s$-compact. Therefore, each assumption of Theorem 1 is satisfied and hence
we have
$$\sup_{x\in X}\inf_{\lambda\in I}(\varphi(T(x)-\lambda x)+\psi(\lambda))=\inf_{\lambda\in I}\sup_{x\in X}
(\varphi(T(x)-\lambda x)+\psi(\lambda))\ .\eqno{(3)}$$
Now, observe that if $\lambda\in I\setminus \{0\}$, we have
$$\sup_{x\in X}\varphi(T(x)-\lambda x)=+\infty\ .\eqno{(4)}$$
Indeed, since the $\tau_s$-interior of $X$ is non-empty and $E$ is reflexive and infinite-dimensional, $X$ turns out to be unbounded. But $T(X)$ is
bounded (since $\varphi$ is coercive) and so, since $\lambda\neq 0$,
$$\sup_{x\in X}\|T(x)-\lambda x\|=+\infty$$
which yields $(4)$ by the coercivity of $\varphi$ again. At this point, the conclusion follows directly from $(2)$, $(3)$ and $(4)$.\hfill
$\bigtriangleup$
\medskip
REMARK 2. - Notice that both infinite-dimensionality of $E$ and compactness of $T$ cannot be dropped in Theorem 2. In this connection, it is enough
to take $T(x)=x$, $I=[0,1]$, $\varphi(x)=\|x\|$, and $\psi=0$.\par
\medskip
REMARK 3. - At present, we do not know any example showing that the reflexivity of $E$ cannot be dropped. However, we conjecture that such an example
can be constructed in infinite-dimensional Banach spaces with the Schur property.
\bigskip
\bigskip
{\bf Acknowledgement.} The author has been supported by the Gruppo Nazionale per l'Analisi Matematica, la Probabilit\`a e le loro Applicazioni (GNAMPA) 
of the Istituto Nazionale di Alta Matematica (INdAM).\par

\vfill\eject
\centerline {\bf References}\par
\bigskip
\bigskip
\noindent
[1]\hskip 5pt B. RICCERI, {\it Some topological mini-max theorems via
an alternative principle for multifunctions}, Arch. Math. (Basel),
{\bf 60} (1993), 367-377.\par
\smallskip
\noindent
[2]\hskip 5pt B. RICCERI, {\it Applications of a theorem concerning sets
with connected sections}, Topol. Methods Nonlinear Anal.,
{\bf 5} (1995), 237-248.\par
\smallskip
\noindent
[3]\hskip 5pt B. RICCERI, 
{\it Nonlinear eigenvalue problems},  
in ``Handbook of Nonconvex Analysis and Applications'' 
D. Y. Gao and D. Motreanu eds., 543-595, International Press, 2010.\par
\smallskip
\noindent
[4]\hskip 5pt S. SIMONS, {\it Minimax theorems and their proofs},
in ``Minimax and Applications'', D.-Z. Du and P. M. Pardalos eds., 1-23,
Kluwer Academic Publishers, 1995.\par

\bigskip
\bigskip
\bigskip
\bigskip
Department of Mathematics\par
University of Catania\par
Viale A. Doria 6\par
95125 Catania\par
Italy\par
{\it e-mail address}: ricceri@dmi.unict.it

\bye